\documentclass[12pt]{article}
\usepackage{amssymb,amsmath}
\usepackage{latexsym, amsfonts, amsmath}

\newcommand{\qed}{\hskip 5mm \rule{2.5mm}{2.5mm}\vskip 10pt}
\newcommand{\R}{{\mathbb R}}

\newcommand{\N}{{\mathbb N}}
\newcommand{\Z}{{\mathbb Z}}
\newcommand{\E}{{\mathbb E}}

\newcommand{\proof}{\noindent{\em Proof:\ }}

\begin{document}
\newtheorem{thm}{Theorem}[section]
\newtheorem{defs}[thm]{Definition}
\newtheorem{lem}[thm]{Lemma}
\newtheorem{note}[thm]{Note}
\newtheorem{cor}[thm]{Corollary}
\newtheorem{prop}[thm]{Proposition}
\renewcommand{\theequation}{\arabic{section}.\arabic{equation}}
\newcommand{\newsection}[1]{\setcounter{equation}{0} \section{#1}}
\title{Mixingales on Riesz spaces}
\author{
 Wen-Chi Kuo\footnote{Funded by an FRC short term postdoctoral fellowship and by NRF Scarce Skills Postdoctoral Fellowship 82395}\\
School of Computational and Applied Mathematics\\
 University of the Witwatersrand\\
 Private Bag 3, P O WITS 2050, South Africa \\ \\
 Jessica Vardy\footnote{Funded in part by the FRC} \&
 Bruce A. Watson\footnote{Funded in part by NRF grant IFR2011032400120 and the Centre for Applicable 
Analysis and Number Theory} \\ 
 School of Mathematics\\
 University of the Witwatersrand\\
 Private Bag 3, P O WITS 2050, South Africa }
\maketitle
\abstract{ \noindent
 A mixingale is a stochastic process which combines properties of martingales and mixing sequences.
 McLeish introduced the term mixingale at the $4^{th}$ Conference on Stochastic Processes and
 Application, at York University, Toronto, 1974, in the context of $L^2$.
In this paper we generalize the concept of a mixingale to the measure-free Riesz space setting
(this generalizes all of the $L^p, 1\le p\le \infty$ variants) and prove that a weak law of large
 numbers holds for Riesz space mixingales. In the process we also generalize the concept of 
uniform integrability to the Riesz space setting.
}
\parindent=0cm
\parskip=0.5cm
\newsection{Introduction}

Mixingales were first introduced by D.L. McLeish in \cite{mcleish1}.  
Mixingales are a generalization of martingales and mixing sequences.  McLeish 
defines mixingales using the $L^2$-norm.  In \cite{mcleish1} McLeish proves
 invariance principles under strong mixing conditions.  In \cite{mcleish} 
a strong law for large numbers is given using mixingales with restrictions 
on the mixingale numbers.

In 1988, Donald W. K. Andrews used mixingales to present $L^1$ and weak laws 
of large numbers, \cite{andrews}.  Andrews used an analogue of McLeish's mixingale 
condition to define $L^1$-mixingales.  The $L^1$-mixingale condition is 
weaker than McLeish's mixingale condition.  Furthermore, Andrews makes no restriction 
on the decay rate of the mixingale numbers, as was assumed by McLeish.  The proofs presented 
in Andrews are remarkably simple and self-contained. Mixingales have also been considered in a
general $L^p, 1\le p<\infty$, by, amongst others de Jong, in \cite{dejong-1, dejong-2}
and more recently by Hu, see \cite{hu}.

In this paper we define mixingales in a Riesz space and prove a weak law of large numbers for
mixingales in this setting.  This generalizes the results in the $L^p$ setting to a measure free
setting.  In our approach the proofs rely on the order structure of the Riesz spaces which
highlights the underlying mechanisms of the theory. This develops on the work of Kuo, Labuschagne,
Vardy and Watson, see \cite{ klw-indag, klw-exp, vw-1, vw-2}, in formulating the theory of stochastic
processes in Riesz spaces.  Other closely related generalizations were given by Stoica \cite{stoica-4}
and Troitsky \cite{troitsky}.

In Section 2 we give a summary of the Riesz space concepts needed as well as the essentials of the
formulation of stochastic processes in Riesz spaces.  Analogous concepts in the classical probability
setting can be found in \cite{billingsley}.  Mixingales in Riesz spaces are defined in Section 3 and 
some of their basic properties derived.  The main result, the weak law of large number for 
mixingales, Theorem \ref{WLLN},  is proved in Section 4 along with a result on the Ces\`aro
summability of martingale difference sequences.

We thank the referees for their valuable suggestions.

\newsection{Riesz Space Preliminaries}

In this section we present the essential aspects of Riesz spaces required for this paper,
 further details can be found in \cite{A-B} or \cite{zaanen}.
The foundations of stochastic processes in Riesz spaces will also given.  
These will form the framework in which mixingales will
be studied in later sections.

A Riesz space, $E$, is a vector space over $\R$, with an order structure that is compatible 
with the algebraic structure on it, i.e. if $f,g\in E$ with $f\le g$ then $f+h\le g+h$ and
$\alpha f\le \alpha g$ for $h\in E$ and $\alpha \ge 0, \alpha \in \R$.
A Riesz space, $E$,  is Dedekind complete if every non-empty upwards directed subset of $E$
 which is bounded above has a supremum.
A Riesz space, $E$, is Archimedean if for each  
$u \in E_+ := \{ f\in E | f\geq 0\}$, the positive cone of $E$, the sequence
$(nu)_{n\in\N}$ is bounded if and only if $u=0$.
We note that every Dedekind complete Riesz space is Archimedean, \cite[page 63]{zaanen}.

We recall from \cite[page 323]{A-B}, for the convenience of the reader, the definition of order
convergence of an order bounded net $(f_\alpha)_{\alpha\in\Lambda}$ 
in a Dedekind complete Riesz space:
$(f_\alpha)$ is order convergent if and only if 
$$\lim\sup_{\alpha} f_\alpha=\lim\inf_{\alpha} f_\alpha.$$
Here 
\begin{eqnarray*}
 \lim\sup_{\alpha} f_\alpha&=&\inf\{ \sup\{ f_\alpha\,|\, \alpha\ge \beta\}\,|\,\beta\in\Lambda\},\\
 \lim\inf_{\alpha} f_\alpha&=&\sup\{ \inf\{ f_\alpha\,|\, \alpha\ge \beta\}\,|\,\beta\in\Lambda\}.
\end{eqnarray*}

Bands and band projections are fundamental to the methods used in our study. 
A non-empty linear subspace $B$ of $E$ is a band if the following conditions are satisfied:
\begin{description}
  \item[(i)]
     the order interval $[-|f|, |f|]$ is in $B$ for each $f\in B$;
  \item[(ii)]
     for each $D\subset B$ with $\sup D\in E$ we have $\sup D \in B$.
\end{description}
The above definition is equivalent to saying that a band is a solid order closed 
vector subspace of $E$.
The band generated by a non-empty subset $D$ of $E$ is the 
intersection of all bands of $E$ containing $D$, i.e. the minimal band containing $D$.
A principal band is a band generated by a 
single element. If $e\in E_+$ and the band generated by $e$ is $E$, then $e$ is called a weak order
unit of $E$ and we denote the space of $e$ bounded elements of $E$ by
$$E^e=\{ f\in E\,:\, |f|\le ke \mbox{ for some } k\in \R_+\}.$$ 
 In a Dedekind complete Riesz space with weak order unit every band is a principal band
 and, for each band $B$  and $u \in E^+$, 
	\[P_Bu: = \sup\{v\,:\,  0\leq v \leq u, v\in B\}\]
exists. 
The above map $P_B$ can be extended to $E$ by setting $P_Bu=P_Bu^+-P_Bu^-$ for $u\in E$.
With this extension, $P_B$ is a positive linear projection which commutes with the operations
of supremum and infimum in that $P(u\vee v)=Pu\vee Pv$ and  $P(u\wedge v)=Pu\wedge Pv$.
Moreover $0\le P_Bu\le u$ for all $u\in E_+$ and the range of $P_B$ is $B$.

In order to study stochastic processes on Riesz spaces, we need to recall the definition of
 a conditional expectation operator on a Riesz space from \cite{klw-indag}.  
As only linear operators between Riesz spaces will be considered, we use the term operator
to denote a linear operator between Riesz spaces.  Let $T: E \to F$ be an operator 
where $E$ and $F$ are Riesz spaces.  We say that $T$ is a positive operator if $T$ maps the 
positive cone of $E$ to the positive cone of $F$, denoted $T\ge 0$.

In this paper we are concerned with order continuous positive operators between Riesz spaces.

\begin{defs}
  Let $E$ and $F$ be Riesz spaces and $T$ be a positive operator between $E$ and $F$.
  We say that $T$ is order continuous if for each directed set
  $D \subset E$ with $f \downarrow_{f\in D}  0$ in $E$ we have that $Tf \downarrow_{f\in D} 0$.

 Here a set $D$ in $E$ is said to be downwards directed if for $f,g\in D$ there exists $h\in D$
 with $h\le f \wedge g$.  In this case we write $D\downarrow$ or $f \downarrow_{f\in D}$. 
 If, in addition, $g=\inf D$ in $E$,
 we write $D\downarrow g$ or $f \downarrow_{f\in D} g$.
 \end{defs}

Note that if $T$ is a positive order continuous operator with $0\leq S \leq T$ (i.e.
 $0\le Sg\le Tg$ for all $g\in E$)
then $S$ is order continuous.   In particular band projections are order continuous. 

\begin{defs}
 Let $E$ be a Dedekind complete Riesz space with weak order unit, $e$.
 We say that $T$ is a conditional expectation operator in $E$ if
 $T$ is a postive order continuous projection which maps weak order units to weak order units
 and has range, ${\cal R}(T)$, a Dedekind complete Riesz subspace of $E$.
\end{defs}

If $T$ is a conditional expectation operator on $E$, as $T$ is a projection it is easy to verify that
at least one of the weak order units of $E$ is invariant under $T$. Various authors have studied
stochastic processes and conditional expectation type operators in terms of order (i.e. in Riesz 
spaces and Banach lattices), see for example \cite{LdP, stoica-4, troitsky}.

To access the averaging properties of conditional expectation operators a multiplicative structure is
needed.  In the Riesz space setting the most natural multiplicative structure is that of
an $f$-algebra. This gives a multiplicative structure that is compatible with 
the order and additive structures on the space.  The space $E^e$, where $e$ is a weak order
unit of $E$ and $E$ is Dedekind complete, has a natural $f$-algebra structure generated by setting
$(Pe)\cdot (Qe)=PQe=(Qe)\cdot (Pe)$ for band projections $P$ and $Q$.
Now using Freudenthal's Theorem this multilpication can be extended to the whole of $E^e$
and in fact to the universal completion $E^u$. Here $e$ becomes the multiplicative unit.
This multiplication is associative, 
distributive and is positive in the sense that if $x,y\in E_+$ then $xy\ge 0$.

If $T$ is a conditional expectation operator on the Dedekind complete Riesz space $E$ with weak 
order unit $e  = Te$, then restricting our attention to the $f$-algebra $E^e$
$T$ is an averaging operator, i.e. $T(fg)=fTg$ for $f,g\in E^e$ and $f\in R(T).$
In fact $E$ is an $E^e$ module which allows the extension of the averaging 
property, above, to $f,g\in E$ with at least one of them in $E^e$.
For more information about $f$-algebras and the averaging properties of conditional 
expecation operators we refer the reader to \cite{BBT, BvR, dodds, grobler, klw-exp}.

\begin{defs}
 Let $E$ be a Dedekind complete Riesz space with weak order unit and $T$ be a strictly positive conditional expectation on $E$.
 The space $E$ is universally complete with respect to $T$, i.e. $T$-universally complete,
 if for each increasing net $(f_\alpha)$ in $E_+$
 with $(Tf_\alpha)$ order bounded, we have that $(f_\alpha)$ is order convergent.
\end{defs}

If $E$ is a Dedekind complete Riesz space and $T$ is a strictly positive conditional 
expectation operator on $E$, then $E$ 
has a $T$-universal completion, see \cite{klw-exp}, which is 
the natural domain of $T$, denoted ${\rm dom}(T)$ in the universal completion, $E^u$, of $E$,
 also see \cite{dodds, grobler}.
Here ${\rm dom}(T)=D-D$ and $Tx:=Tx^+-Tx^-$ for $x\in {\rm dom}(T)$ where
$$D=\{ x\in E^u_+ | \exists(x_\alpha)\subset E_+, x_\alpha\uparrow x, (Tx_\alpha)
   \ \mbox{order bounded in}\ E^u\},$$
and $Tx:=\sup_\alpha Tx_\alpha$, for $x\in D$, where $(x_\alpha)$ is an increasing net in $E_+$
 with $(x_\alpha)\subset E_+$, $(Tx_\alpha)$ order bounded in $E^u$.

Building on the concept of a conditional expectation, if $(T_i)$ is a sequence of conditional
expectations on $E$ indexed by either $\N$ or $\Z$, we say that $(T_i)$
is a filtration on $E$ if
$$T_iT_j=T_i=T_jT_i,\quad \mbox{for all}\quad i\le j.$$
If $(T_i)$ is a filtration and $T$ is a conditional expectation with $T_iT=T=TT_i$ for all $i$, 
then we say that the filtration is compatible with $T$.
The sequence $(T_i)$ of conditional expectations in $E$ being a filtration 
is equivalent to ${\cal R}(T_i)\subset {\cal R}(T_j)$ for $i\le j$.
If $(T_i)$ is a filtration on $E$ and $(f_i)$ is a sequence in $E$, we say that 
$(f_i)$ is adapted to the filtration $(T_i)$ if $f_i\in{\cal R}(T_i)$ for all $i$ in the
index set of the sequence $(f_i)$.  The double sequence $(f_i,T_i)$ is called a martingale
if $(f_i)$ is adapted to the filtration $(T_i)$ and in addition
$$f_i=T_if_j,\quad\mbox{for}\quad i\le j.$$
The double sequence $(g_i,T_i)$ is called a martingale difference sequence
if $(g_i)$ is adapted to the filtration $(T_i)$ and
$$T_i g_{i+1}=0.$$
We observe that if $(f_i)$ is adapted to the filtration $(T_i)$ then $(f_i-T_{i-1}f_i,T_i)$
is a martingale difference sequence.  Conversely, if $(g_i,T_i)$ is a martingale difference
sequence, then $(s_n,T_n)$ is a martingale, where
$$s_n=\sum_{i=1}^n g_i, \quad n\ge 1,$$
and the martingale difference sequence generated from $(s_n,T_n)$ is
precisely $(g_n,T_n)$.

We now give some basic aspects of independence in Riesz spaces.  
An in depth discussion of independence in the context of Markov processes in
Riesz spaces can be found in \cite{vw-1, vw-2}.

\begin{defs}
 Let $E$ be a Dedekind complete Riesz space with conditional expectation $T$ and weak order unit $e=Te$.
 Let $P$ and $Q$ be band projections on $E$, we say that $P$ and $Q$ are $T$-conditionally independent with respect to $T$ if
 \begin{eqnarray}
  TPTQe=TPQe=TQTPe.\label{indep-e}
 \end{eqnarray}
 We say that two Riesz subspaces $E_1$ and $E_2$ of $E$ containing ${\cal R}(T)$,
 are $T$-conditionally independent with respect to $T$ if
 all band projections $P_i, i=1,2,$ in $E$ with  $P_ie\in E_i, i=1,2,$ 
 are $T$-conditionally independent with respect to $T$.
\end{defs}

\begin{cor}
 Let $E$ be a Dedekind complete Riesz space with conditional expectation $T$ and let $e$ 
 be a weak order unit which is invariant under $T$.
 Let $P_i, i=1,2,$ be band projections on $E$. Then $P_i, i=1,2,$  are $T$-conditionally independent
 if and only if the closed Riesz subspaces
 $E_i=\left< P_ie, \mathcal{R}(T)\right>, i=1,2,$ are $T$-conditionally independent.
\end{cor}

If $(\Omega, {\cal A}, \mu)$ is a probability space and $f_\alpha, \alpha\in\Lambda,$ 
is a family in
$L^1(\Omega, {\cal A}, \mu)$, indexed by $\Lambda$, 
the family is said to be uniformly integrable if
for each $\epsilon>0$ there is $c>0$ so that 
$$\int_{\Omega_\alpha (c)} |f_\alpha|\,d\mu\le \epsilon,\quad\mbox{for all}\quad \alpha\in\Lambda,$$
where $\Omega_\alpha (c)=\{ x\in\Omega \,:\, |f_\alpha (x)| > c\}.$
This concept can be extended to the Riesz space setting as
$T$-uniformity, 
see the definition below,
where $T$ is a conditional expectation operator.  
In the case of  the Riesz space being $L^1(\Omega, {\cal A}, \mu)$ and $T$ being the 
expectation operator, the two concepts coincide.

\begin{defs}
Let $E$ be a Dedekind complete Riesz space with conditional expectation operator
$T$ and weak order unit $e=Te$. 
 Let $f_\alpha,\alpha\in\Lambda,$ be a family in $E$, where $\Lambda$ is some index set.
 We say that $f_\alpha,\alpha\in\Lambda,$ is $T$-uniform if 
\begin{eqnarray}
 \sup\{ TP_{(|f_\alpha|-ce)^+}|f_\alpha| \,:\, \alpha\in\Lambda\} \to 0 \quad \text{ as } \quad c\to \infty.
\label{ui}
\end{eqnarray}
\end{defs}

\begin{lem}\label{UI}
 Let $E$ be a Dedekind complete Riesz space with conditional expectation $T$ and let $e$ 
 be a weak order unit which is invariant under $T$.
 If $f_\alpha\in E, \alpha\in\Lambda,$ is a $T$-uniform family, then 
 the set $\{ T|f_\alpha|\,:\, \alpha\in\Lambda\}$ is bounded in $E$.
\end{lem}

\proof
 As the sequence $f_\alpha,\alpha\in\Lambda,$ is $T$-uniform
 \begin{eqnarray*}
  J_c:=\sup\{ TP_{(|f_\alpha|-ce)^+}|f_i|\, :\, \alpha\in\Lambda\}
   \to 0 \quad \text{ as } \quad c\to \infty.
 \end{eqnarray*}
 In particular this implies that $J_c$ exists for $c>0$ large and that, for sufficiently large $K>0$,
 the set $\{J_c\,:\, c\ge K\}$ is bounded in $E$.
 Hence there is $g\in E_+$ so that
 \begin{eqnarray*}
 TP_{(|f_\alpha|-ce)^+}|f_\alpha|\le g,\quad\mbox{for all}\quad \alpha\in\Lambda, c\ge K,
 \end{eqnarray*}
 By the definition of  $P_{(|f_\alpha| -ce)^+}$,
 $$(I-P_{(|f_\alpha| - ce)^+})|f_\alpha|\le ce,\quad\mbox{for} \quad \alpha\in\Lambda, c>0.$$
 Combining the above for $c=K$ gives
\begin{eqnarray*}
  T|f_\alpha| =TP_{(|f_\alpha| - Ke)^+}|f_\alpha| + T(I-P_{(|f_\alpha| - Ke)^+})|f_\alpha|
  \leq g + Ke,
\end{eqnarray*}
for all $\alpha\in\Lambda$.
\qed

\newsection{Mixingales in Riesz Spaces}

In terms of classical probability theory  $((f_i)_{i\in\N},({\cal A}_i)_{i\in \Z})$ is a mixingale in the 
probability space $(\Omega,{\cal A},\mu)$ if $({\cal A}_i)_{i\in \Z}$ is an increasing sequence of 
sub-$\sigma$-algebras of ${\cal A}$ and $(f_i)_{i\in \N}$ is a sequence of
 ${\cal A}$ measurable functions with
$\E[|\E[f_i|{\cal A}_{i-m}]|]$ and $\E[|f_i-\E[f_i|{\cal A}_{i+m}]|]$ existing and 
$\E[|\E[f_i|{\cal A}_{i-m}]|]\le c_i\Phi_m$ and $\E[|f_i-\E[f_i|{\cal A}_{i+m}]|]\le c_i\Phi_{m+1}$
for some sequences $(c_i), (\Phi_i)\subset \R_+$ with $\Phi_i\to 0$ as $i\to\infty$.
They were first introduced, in \cite{mcleish}, with the additional assumption that 
$(f_i)\subset L^2(\Omega,{\cal A},\mu)$, this was later generalized, in \cite{andrews}, to
$(f_i)\subset L^1(\Omega,{\cal A},\mu)$.  We now formulate a measure free abstract
definition of a mixingale in the setting of Riesz spaces with conditional expectation operator.
This generalizes the above classical definitions.

\begin{defs}\label{mixing}
Let $E$ be a Dedekind complete Riesz space with conditional expectation operator, $T$, 
and weak order unit $e = Te$.  
Let $(T_i)_{i\in \Z}$ be a filtration on $E$ compatible
with $T$ in that $T_iT=T=TT_i$ for all $i\in\Z$.
Let $(f_i)_{i\in\N}$ be a sequence in $E$.  
We say that $(f_i, T_i)$ is a 
mixingale in $E$ compatible with $T$
if there exist  $(c_i)_{i\in\N}\subset E_+$ and $(\Phi_m)_{m\in\N}\subset \R_+$ 
such that $\Phi_m \to 0$ as $m \to \infty$ and for all $i,m\in\N$ we have
	\begin{itemize}
	\item[(i)] $T|T_{i-m}f_i|\leq\Phi_m c_i$,
	\item[(ii)] $T|f_i - T_{i+m}f_i| \leq \Phi_{m+1}c_i$.
	\end{itemize}
\end{defs}

The numbers $\Phi_m, m\in\N,$ are referred to as the mixingale numbers.
These numbers give a measure of the temporal dependence of the sequence $(f_i)$.
The constants $(c_i)$ are chosen to index the \lq magnitude\rq \  of the the 
random variables $(f_i)$. 

In many applications the sequence $(f_i)$ is adapted to the filtration $(T_i)$.
The following theorem sheds more light on the structure of mixingales for this special case.

 We recall that if $T$ is a conditional expectation operator on a Riesz space $E$
 then $T|g|\ge |Tg|$. 

\begin{lem}
 Let $E$ be a Dedekind complete Riesz space with 
conditional expectation operator, $T$, 
and weak order unit $e = Te$.
Let $(f_i, T_i)_{i\geq 1}$ be a mixingale in $E$ compatible with $T$.
\begin{itemize}
  \item[(a)]{The sequence $(f_i)$ has $T$-mean zero, i.e. $Tf_i=0$ for all $i\in\N$.}
  \item[(b)]{If in addition $(f_i)_{i\in\N}$ is $T$-conditionally independent and
                ${\cal R}(T_i) = \left< f_1, \dots, f_{i-1}, {\cal R}(T)\right>$ 
                  then the mixingale numbers may be taken as zero, where
                 $\left<f_1,\dots,f_{i-1},{\cal R}(T)\right>$ is the order closed Riesz subspace of $E$ generated
                 by $f_1,\dots,f_{i-1}$ and ${\cal R}(T)$.}
\end{itemize}
\end{lem}

\proof
{\bf (a)}  Here we observe that the index set for the filtration $(T_i)$ is $\Z$, thus 
 \begin{eqnarray*}
  |Tf_i| &=& |TT_{i-m}f_i|\\
   &\leq & T|T_{i-m} f_i|\\
   &\leq & \Phi_m c_i\\
   &\to& 0 \quad \text{ as } m \to \infty
 \end{eqnarray*}
 giving $Tf_i=0$ for all $i\ge 0$.

{\bf (b)} As $f_i \in {\cal R}(T_{i+1}) \subset {\cal R}(T_{i+m})$, it follows that 
  $$f_i- T_{i+m}f_i = 0,\quad\mbox{ for all }\quad i,m\in\N.$$
We recall from \cite{vw-1} that two closed Riesz spaces, say $E_1$ and $E_2$, are 
$T$-conditionally independent if and only if $T_if = Tf$ for all $f \in E_{3-i},\  i =1,2$.   As $(f_i)$ is $T$-conditionally 
independent and
 as $(f_i)$ has $T$-mean zero (from (a)), we have that
 $$T_{i-m}f_i=Tf_i=0,$$
 for $i,m\in\N$
 Thus we can choose $\Phi_m = 0$ for all $m \in\N$.  
\qed

\newsection{The Weak Law of  Large Numbers}

We now show that the above generalization of mixingales to the measure free Riesz
space setting admits a weak law of large numbers.  For a sequence 
$(f_i)$, we shall denote its Cez\`{a}ro mean by 
	\[\overline{f}_n = \frac1n\sum_{i=1}^n f_i.\]

\begin{lem}\label{the-lemma}
Let $E$ be a Dedekind complete Riesz space with conditional expectation operator $T$,
weak order unit $e= Te$ and filtration $(T_i)_{i\in \N}$ compatible with $T$.  
Let $(f_i)$ be an $e$-uniformly bounded sequence adapted to the filtration $(T_i)$,
and $g_i:=f_i-T_{i-1}f_i$,
 then $(g_i,T_i)$ is a martingale difference sequence with 
	\[T|\overline{g}_n| \to 0 \quad \text{ as } n\to \infty.\]
\end{lem}

\proof
 Clearly 
 $$T_ig_{i+1}=T_if_{i+1}-T^2_if_{i+1}=0$$
 and $(g_i)$ is adapted to $(T_i)$ so indeed $(g_i,T_i)$ is a martingale difference sequence.
 
 Let $B>0$ be such that $|f_i|\le Be$, for all $i\in\N$.
 For $j>i$, as $T_jT_i=T_i$ and $T_if_i=f_i$ it follows that $T_ig_i=g_i$ and
 $$T_{i}g_j=T_{i}f_j-T_{i}T_{j-1}f_j=T_{i}f_j-T_{i}f_j=0.$$
 In addition, 
 $$|g_i|\le |f_i|+|T_{i-1}f_i|\le |f_i|+T_{i-1}|f_i|\le 2Be,\quad\mbox{for all}\quad i\in\N.$$

 Set 
 $$s_n=\sum_{i=1}^n g_i.$$
 As $|g_i|\le 2Be$ we have that $g_i$ is in the $f$-algebra $E^e$.  Hence the product
 $g_ig_j$ is defined in $E^e$. 
 Now as $T_j$ is an averaging operator, see \cite{klw-exp}, and $g_j\in{\cal R}(T_j)$ we have
 $$T_i(g_ig_j)=g_iT_ig_j=0,\quad \mbox{for}\quad j>i.$$

 Combining these results gives
 \begin{eqnarray*}
  T(s_n^2)&=&\sum_{i,j=1}^n T(g_ig_j)\\
     &=& \sum_{i=1}^n T(g_i^2)+2\sum_{i<j} T(g_ig_j)\\ 
     &=& \sum_{i=1}^n T(g_i^2)+2\sum_{i<j} TT_i(g_ig_j)\\ 
     &=& \sum_{i=1}^n T(g_i^2)+0. 
 \end{eqnarray*}

 Thus 
 $$T(s_n^2)=\sum_{i=1}^n T(g_i^2).$$
 But 
 $$g^2_i=|g_i|^2\le 4B^2e$$
 as $e$ is the algebraic identity of the $f$-algebra $E^e$ and $|g_i|\le 2Be$. Thus
 \begin{eqnarray}
  T(s_n^2)\le 4nB^2e.\label{square}
 \end{eqnarray}
 Now let
 $$J_n=P_{s_n^+}-(I-P_{s_n^+})$$
 where $P_{s_n^+}$ is the band projection on the band in $E$ generated by $s_n^+$.
 From the definition of the $f$-algebra structure on $E^e$, if $P$ and $Q$ are band projections 
 then $(Pe)(Qe)=PQe$ which together with Freudenthal's Theorem \cite[page 216]{zaanen} 
 enables us to conclude $$|s_n|=J_ns_n=(J_ne)s_n$$
 and $(J_ne)^2=e$, as $J_n^2=I$.  But
\begin{eqnarray*}
 0&\le&\left(\frac{J_ne}{n^{1/4}}-\frac{s_n}{n^{3/4}}\right)^2\\
   &=&\left(\frac{J_ne}{n^{1/4}}\right)^2 + 
          \left(\frac{s_n}{n^{3/4}}\right)^2
          -2\frac{J_ne}{n^{1/4}}\frac{s_n}{n^{3/4}}\\
   &=&\frac{e}{n^{1/2}} + 
          \frac{s_n^2}{n^{3/2}}
          -2\frac{|s_n|}{n}
 \end{eqnarray*}
 giving
 $$\frac{e}{n^{1/2}} + \frac{s_n^2}{n^{3/2}}\ge 2\frac{|s_n|}{n}.$$

 Applying $T$ to this inequality gives
 $$\frac{e}{n^{1/2}} + 
       \frac{T(s_n^2)}{n^{3/2}}
          \ge 2\frac{T|s_n|}{n}$$
 Combining the above inequality with (\ref{square}) gives
 \begin{eqnarray*}
  2\frac{T|s_n|}{n}&\le&\frac{e}{n^{1/2}} +  \frac{T(s_n^2)}{n^{3/2}}\\
          &\le&\frac{e}{n^{1/2}} + \frac{4nB^2e}{n^{3/2}}\\
          &=&\frac{1 + 4B^2}{n^{1/2}}  e,
 \end{eqnarray*}
 and thus
 \begin{eqnarray}
  T|\overline{g}_n|
          &\le&\frac{1 + 4B^2}{2n^{1/2}}  e.\label{g-bound}
 \end{eqnarray}
 Since $E$ is an Archimedean Riesz space letting $n\to\infty$ in (\ref{g-bound}) gives
  $T|\overline{g}_n|\to 0$ as $n\to\infty$.
\qed


We are now able to prove an analogue to the weak law of large numbers for
mixingales in Riesz spaces.

\begin{thm}\label{WLLN}{\rm\bf [Weak Law of Large Numbers]}
Let $E$ be a Dedekind complete Riesz space with conditional expectation operator $T$,
weak order unit $e= Te$ and filtration $(T_i)_{i\in \Z}$.  Let $(f_i, T_i)_{i\geq 1}$ be a 
$T$-uniform mixingale with $c_i$ and $\Phi_i$ as defined in Definition \ref{mixing}.
\begin{itemize}
\item[(a)]
If $\displaystyle{\left( \frac{1}{n} \sum_{i=1}^n c_i \right)_{n\in \N}}$ is bounded in $E$ then 
 \[ T|\overline{f}_n| = T\left|\frac{1}{n}\sum_{i=1}^n f_i\right| \to 0 \quad \text{ as } n \to \infty.\]
\item[(b)]
If $c_i = T|f_i|$ for each $i\geq 1$ then 
 \[ T|\overline{f}_n| = T\left|\frac{1}{n}\sum_{i=1}^n f_i\right| \to 0 \quad \text{ as } n \to \infty.\]
\end{itemize}
\end{thm}

\proof	
 {\bf (a)}
  Let $$y_{m,i} = T_{i+m}f_i - T_{i+m - 1}f_i,\quad\mbox{for}\quad i \geq 1, m \in \Z.$$
  Let $B> 0$.  Let $h_i=(I-P_{(|f_i|-Be)^+})f_i$ and $d_i=P_{(|f_i|-Be)^+}f_i, i\in \N$ then
  $f_i=h_i+d_i$.
  Now $(T_{i+m}h_i)_{i\in\N}$ is $e$-bounded and adapted to $(T_{i+m})_{i\in\N}$, so 
  from Lemma \ref{the-lemma}
  $(T_{i+m}h_i-T_{i+m-1}h_i,T_{i+m})_{i\in\N}$  
  is a martingale difference sequence with
  $$T\left|\frac{1}{n}\sum_{i=1}^n (T_{i+m}h_i-T_{i+m-1}h_i)\right|\to 0$$  
  as $n\to\infty$.
  Using the $T$-uniformity of $(f_i)$ we have
  \begin{eqnarray*}
  T\left|\frac{1}{n}\sum_{i=1}^n (T_{i+m}d_i-T_{i+m-1}d_i)\right|  
   &\le& T\frac{1}{n}\sum_{i=1}^n |T_{i+m}d_i-T_{i+m-1}d_i|\\  
   &\le& \frac{1}{n}\sum_{i=1}^n (TT_{i+m}|d_i|+TT_{i+m-1}|d_i|)\\  
   &=& \frac{2}{n}\sum_{i=1}^n T|d_i|\\  
   &\le& 2\sup\{ T|d_i|\,:\, i=1,\dots,n\}.
  \end{eqnarray*}  
  Thus
 \begin{eqnarray*}
  T\left|\frac{1}{n}\sum_{i=1}^n (T_{i+m}d_i-T_{i+m-1}d_i)\right|  
   &\le& 2\sup\{ T|P_{(|f_i|-Be)^+}f_i|\,:\, i\in\N\}.
  \end{eqnarray*}  
  Combining the above results gives 
  \begin{eqnarray*}
  \limsup_{n\to\infty}T\left|\frac{1}{n}\sum_{i=1}^n (T_{i+m}f_i-T_{i+m-1}f_i)\right|  
   &\le& 2\sup\{ TP_{(|f_i|-Be)^+}|f_i|\,:\,i\in\N\}\\
  &\to&0
  \end{eqnarray*}  
  as $B\to\infty$ by the $T$-uniformity of $(f_i)$.
  Thus $T|\overline{y}_{m,n}| \to 0$ as $ n\to \infty$.

We now make use of a telescoping series to expand $\overline{f}_n$, 
\begin{eqnarray*}
\overline{f}_n &=& \frac1n \sum_{i=1}^n f_i\\
&=& \frac{1}{n} \sum_{i=1}^n\left(f_i - T_{i+M}f_i 
     + \sum_{m = -M +1}^M (T_{i+m}f_i - T_{i+m-1}f_i)
	+ T_{i-M}f_i\right) \\
&=& \frac1n \sum_{i=1}^n(f_i - T_{i+M}f_i) + \sum_{m = -M +1}^M \frac1n\sum_{i=1}^n(T_{i+m}f_i - T_{i+m-1}f_i)
	+ \frac1n\sum_{i=1}^nT_{i-M}f_i\\
&=& \frac1n \sum_{i=1}^n(f_i - T_{i+M}f_i) + \sum_{m = -M +1}^M \overline{y}_{m,n} + \frac1n\sum_{i=1}^nT_{i-M}f_i\\
\end{eqnarray*}
Applying $T$ to the above expression we can bound $T|\overline{f}_n|$ by means of the 
defining properties of a mixingales as follows
\begin{eqnarray*}
T|\overline{f}_n| 
&\leq & \frac1n \sum_{i=1}^nT|f_i - T_{i+M}f_i| + \sum_{m = -M +1}^MT|\overline{y}_{m,n}| 
	+ \frac1n\sum_{i=1}^nT|T_{i-M}f_i|\\
&\leq & \frac1n \sum_{i=1}^n  \Phi_{M+1}c_i + \sum_{m = -M +1}^MT|\overline{y}_{m,n}|+ 
	\frac1n \sum_{i=1}^n \Phi_M c_i.
\end{eqnarray*}
  Since 
$\displaystyle{\left(\frac1n\sum_{i=1}^{n} c_i \right)_{n\in\N}}$ is bounded in $E$ there is 
$q\in E_+$ so that 
$\displaystyle{\frac1n\sum_{i=1}^{n} c_i \le q,}$ for all $n\in\N$, which when combined with the 
above display yields 
\begin{eqnarray*}
T|\overline{f}_n| 
&\leq & (\Phi_{M+1}+\Phi_M)q + \sum_{m = -M +1}^MT|\overline{y}_{m,n}|.
\end{eqnarray*}
Letting $n\to\infty$ gives
\begin{eqnarray*}
\limsup_{n\to\infty}T|\overline{f}_n| 
&\leq & (\Phi_{M+1}+\Phi_M)q.
\end{eqnarray*}
Now taking $M\to\infty$ gives
\begin{eqnarray*}
\limsup_{n\to\infty}T|\overline{f}_n|=0, 
\end{eqnarray*}
completing the proof of (a).

{\bf (b)} By Lemma \ref{UI}, $(T|f_i|)$ is bounded, say by $q\in E_+$, so
\begin{eqnarray*}
\limsup_{n\to\infty} \frac{1}{n}\sum_{i=1}^n c_i 
= \limsup_{n\to\infty}  \frac{1}{n}\sum_{i=1}^n T|f_i|\le q,
\end{eqnarray*}
making (a) applicable.
\qed


\end{document}